\documentclass[10pt]{amsart}
\usepackage{latexsym}
\usepackage{amssymb}
\usepackage{amsmath}
\usepackage{mathrsfs}
\def\ба{\operatorname{ба}}
\def\comp{\operatorname{comp}}

\begin{document}

\title[Sobolev and Schwartz]
{Sobolev and Schwartz:\\
Two Fates and Two Fames}
\author{S.~S. Kutateladze}

\address[]{
Sobolev Institute of Mathematics\newline
\indent 4 Koptyug Avenue\newline
\indent Novosibirsk, 630090
\indent RUSSIA}
\email{
sskut@member.ams.org
}
\dedicatory{
On the Occasion of the Centenary of\\
the Birth of S.~L. Sobolev}
\begin{abstract}
This is a brief overview of the lives and contributions of
S.~L. Sobolev  and L.~Schwartz, the cofounders of
distribution theory.
\end{abstract}
\date{February 4, 2008}
\thanks
{Partly printed in~\cite{Herald} with unauthorized omissions.
\hfill\break
\phantom{qq} The author thanks  V.~A.~Aleksandrov and V.~P.~Golubyatnikov
who helped him in better  understanding of French sources.
The author is especially grateful to Yu.~L.~Ershov who was persistent
in inviting the author to make a talk at the special session
of the Academic Council of the Sobolev Institute on October~14, 2003.  The present article bases
on this talk.
The author acknowledges the subtle and deep comments of V.~I. Arnold and
V.~S. Vladimirov on the preprint of a~draft of the talk which led to
many improvements.
}

\maketitle

In the history of mathematics there are quite a few persons
whom we prefer to recollect in pairs.
Listed among them are Euclid and Diophant, Newton and Leibniz,
Bolyai and Lobachevski\u\i{}, Hilbert and Poincar\'e, as well as
Bourbaki and Arnold. In this series we enroll  Sobolev and Schwartz
who are inseparable from one of the most brilliant discoveries
of the twentieth century, the theory of generalized functions or distribution
theory,  providing a revolutionary new approach to partial differential equations.

The most vibrant and lasting achievements  of mathematics
reside in formulas and lists. There are pivotal distinctions between
lists and formulas. The former deposit that which was open for us.
The lists of platonic solids, elementary catastrophes, and finite
simple groups are next of kin to the {\it Almagest} and herbaria.
They are the objects of admiration, tremendous and awe-struck.
The article of the craft of mathematics is a~formula.
Each formula  enters into life as an instance of materialization of mathematical
creativity. No formula  serves  only  the purpose it was intended to.
In part, any formula is reminiscent of  household appliances,
toys, or software. It is a very rare event that somebody reads
the user's guide of a new TV set or the manual for running a new computer
program. Usually everyone utilizes his or her new gadgets experimentally
by pressing whatever keys and switches.
In much the same way we handle formulas.
We  painstakingly ``twist and turn'' them, audaciously insert new parameters,
willfully interpret  symbols, and so on.

{\sc Mathematics is the craft of formulas and the art
of calculus}. If someone considers this claim as
feeble and incomplete, to remind is in order that, logically speaking,
set theory is just an instance of the first order predicate calculus.

Distribution theory has become the calculus of today.
Of such a scale and scope is the scientific discovery by
Sobolev and Schwartz.

\section{Serge\u\i{} L$'$vovich Sobolev}

Sobolev was born in St. Petersburg on October 6, 1908   in
the family of Lev Aleksandrovich Sobolev, a solicitor.
Sobolev's    grandfather on his father's side descended
from a~family of Siberian Cossacks.

Sobolev was bereaved of his father in the early childhood and was raised
by his mother Natal$'$ya Georgievna who was a highly-educated  teacher
of literature and history. His mother also had the second
speciality: she graduated from a medical institute and worked as a tutor at
the First Leningrad Medical Institute. She cultivated in Sobolev
the decency, indefatigability,  and endurance that
characterized him as a scholar and personality.

Sobolev fulfilled the program of secondary school at home, revealing
his great attraction  to mathematics. During the Civil War he and his mother
lived in Kharkov. When
living there, he studied at the preparatory courses of a evening
technical  school for one semester. At the age of 15 he completed
the obligatory programs of secondary school in  mathematics,
physics, chemistry, and other natural sciences,  read
the classical pieces of the Russian and
world   literature as well as many books on philosophy,
medicine, and biology.

After the family had transferred from Kharkov to Petersburg in 1923, Sobolev
entered the graduate class of
School No.~190 and finished  with honors in 1924, continuing
his study at the First State Art School in the piano class.
At the same year  he entered the Faculty of Physics and Mathematics of
Leningrad State University (LSU) and attended the lectures of Professors
N.~M.~G\"unter,
V.~I.~Smirnov, G.~M.~Fikhtengolts, and others. He made his diploma on the analytic
solutions of a system of differential equations with two independent
variables under the supervision of G\"unter.

G\"unter propounded the idea that the set functions are
inevitable in abstracting the concept of solution to a
differential equation. G\"unter's approach influenced the further train
of thought of
Sobolev.\footnote{It was A.~M. Vershik and V.~I. Arnold who attracted
the author's attention to the especial role of G\"unter
in the prehistory of distribution theory.}

After graduation from LSU in 1929, Sobolev started his work
at the Theoretical Department of the Leningrad
Seismological Institute.
In a close cooperation with Smirnov he then solved some fundamental
problems of wave propagation. It was Smirnov whom Sobolev called
his teacher alongside G\"unter up to his terminal days.

Since 1932 Sobolev  worked at the Steklov Mathematical Institute in
Leningrad; and since 1934, in Moscow. He continued the
study of hyperbolic equations and proposed a new method for
solving the Cauchy problem for a hyperbolic equation with variable
coefficients. This method was based on a generalization
of the Kirchhoff formula.  Research into hyperbolic equations
led Sobolev  to revising  the classical concept of a solution
to a differential   equation. The concept of a generalized or weak
solution of a differential equation    was considered earlier.
However, it was exactly in the works by Sobolev that this
concept was elaborated and applied systematically.
Sobolev posed and solved the Cauchy problem in spaces of functionals,
which was based on the revolutionary extension of the
Eulerian concept of function and declared 1935 as the date of
the birth of the theory of distributions.

Suggesting his
definition of generalized derivative, Sobolev enriched mathematics
with the  spaces of functions whose weak derivatives are integrable
to some   power. These are now called {\it Sobolev spaces}.

\begin{itemize}
\item[]{\small
Let $f$ and $g$ be locally summable functions on
an open subset $G$ of~$\mathbb R^n$, and let $\alpha$ be a multi-index.
A function $g$, denoted by $D^\alpha f$, is the {\it generalized derivative in the Sobolev sense\/}
or {\it weak derivative\/} of $f$ of order $\alpha$
provided that
$$
\int\limits_G f(x)D^\alpha\varphi (x)\,dx=
(-1)^{|\alpha |}
\int\limits_G g(x)\varphi (x)\,dx,
$$
for every {\it test function\/}
$\varphi$, i.~e. such that the support of~$\varphi$
is a compact subset of~$G$ and $\varphi$ is
$|\alpha|=\alpha _1 +\cdots +\alpha _n$ times
continuously differentiable
in ~$G$, where
 $D^\alpha\varphi$ is the classical derivative of
$\varphi$ of order~$\alpha$.
The vector space $W^l_p$, with $p\ge 1$, of the
(cosets of) locally summable $f$ on $G$, whose all weak derivatives
$D^\alpha f$  with $|\alpha |\le l$ are $p$-summable in $G$
becomes a Banach space under the norm:
$$
\|f\|_{W^l_p} =\biggl(\,
\int\limits_G|f|^p\,dx\biggr)^{1/p} +
\sum\limits_{|\alpha |=l}\biggl(\,\int\limits_G
|D^\alpha f|^p\,dx\biggr)^{1/p} .
$$
}
\end{itemize}
\noindent
Sobolev  found the general criteria for equivalence of various norms on
$ W_p^l$ and showed that these spaces are the natural environment for
posing the boundary value problems for elliptic
equations. This conclusion was   based on his thorough study
of the properties of Sobolev spaces. The most important
facts are {\it embedding theorems}. Each embedding
theorem estimates the operator norm of
an embedding, yielding special inequalities between
the norms of one and the same function  inside
various spaces.

The contributions of Sobolev brought him recognition in the USSR.
In 1933 Sobolev was elected a
corresponding member of the Academy of Sciences at the age of 24 years.
In 1939 he became a  full member of the Academy
and  remained the youngest  academician for many years.

Inspired by military applications in the 1940s, Sobolev studying the system
of differential equations describing small oscillations of a rotating fluid.
He obtained   the conditions for stability of a rotating body with
a filled-in cavity which depend on  the shape and parameters of the cavity.
Moreover, he elaborated the cases in which the cavity is
a~ cylinder or an ellipsoid of rotation.
This research by Sobolev  signposted another area of
the general theory which concerns   the Cauchy
and boundary value problems for the equations and systems that are
not solved with respect to  higher time derivatives.

In the grievous years of the Second World War from 1941 to 1944
Sobolev occupied the position of the director of the
Steklov Mathematical Institute.

Sobolev was one of the first scientists who foresaw
the future of   computational mathematics and cybernetics.
From 1952 to   1960 he held the chair of the first national department of
computational mathematics at  Moscow   State University. This
department  has played  a key role in the development
of many important areas of the today's  mathematics.

Addressing the problems of computational mathematics,
Sobolev lavishly  applied the apparatus of the modern sections of the
theoretical core of mathematics. It is typical for him to pose the
problems of     computational mathematics within  functional analysis.
Winged are his words that
``to  conceive the theory of computations without Banach spaces
is impossible just  as trying to conceive it without computers.''

It is worthwhile to emphasize the great role in the
uprise of cybernetics and other new areas of research in this country
which was played  by the publications and speeches of Sobolev who
valiantly defended the new trends in science from the ideologized obscurantism.

To overrate is difficult  the contribution of
Sobolev to the design of the nuclear shield of this
country. From the first stages of the
atomic project of the USSR he was listed among the
top officials of Laboratory No.~2 which
was renamed for secrecy reasons into the Laboratory
of Measuring Instruments (abbreviated as LIPAN in Russian).
Now LIPAN lives as the Kurchatov Center.
The main task of the joint work with I.~K. Kikoin was
the implementation of gaseous diffusive uranium
enrichment for creation of a~nuclear explosive device.

Sobolev administered and supervised various computational teams,
studied the control of the industrial processes of isotope separation,
struggled for the low costs of production and made decisions on many
managerial and technological matters.
For his contribution to the A-bomb project,
Sobolev twice gained a Stalin Prize of the First Degree.
In January of 1952 Sobolev was awarded with  the highest title of
the USSR: he was declared the Hero of the Socialist Labor
for  exceptional service to the state.

Sobolev's research  was inseparable from his management in science.
At the end of the 1950s  M.~A.~Lavrent$'$ev, S.~L.~Sobolev, and
S.~A.~Khristianovich came out with the initiative to organize a new big
scientific center, the Siberian Division of the Academy of Sciences.
For many  scientists of the first enrolment to the Siberian Division
it was  the example of Sobolev, his authority in science,
and the attraction of his personality that yielded
the final argument in deciding to move to Novosibirsk.

The Siberian period of Sobolev's life in science was marked by
the great achievements in the theory of cubature formulas.
Approximate integration is one of the main problems in the
theory of computations---the cost of computation of multidimensional integrals
is extremely high.  Optimizing the formulas of integration
is understood now to be the problem of  minimizing the norm of
the error  on some function space. Sobolev suggested
new approaches to the problem and discovered marvelous classes of
optimal cubature formulas.

Sobolev merits brought him many decorations  and signs of distinction.
In 1988 he was awarded the highest prize of the
Russian Academy of Sciences, the Lomonosov Gold Medal.

Sobolev passed away in Moscow on January 3, 1989.

\section{Laurent Schwartz}

Schwartz was born in Paris on March 5, 1915 in the family of
Anselme Schwartz, \hbox{a~surgeon}. There were quite a few prominent persons among his next of kin.
J.~Hadamard was his granduncle.
Many celebrities are listed in the lone of his  mother's line
Claire Debr\'es: Several  Gaullist politicians belonged to
the Deb\'ers.
In 1938 Schwartz married Marie-H\'el\`en L\'evy, the daughter of
the outstanding mathematician P.~L\'evy who was one of the
forefathers of functional analysis.
Marie-H\'el\`en had become  a professional mathematician and
gained the position of a full professor in~1963.

The munificent gift of Schwartz was revealed in his lec\'ee years.
He won the  most prestigious competition for high school students,
Concours G\'en\'eral in Latin. Schwartz was unsure about his future
career, hovering between geometry and ``classics''  (Greek and Latin).
It is curious that Hadamard had a low opinion of
Schwartz mathematical plans, since the sixteen-years old Laurent
did not know the Riemann zeta function.  By a startling contrast,
~Schwartz was boosted to geometry by the pediatrician Robert Debr\'e
and one of his teachers of classics.


In 1934 Schwartz passed examinations to the \'Ecole Normale Sup\'erieure (ENS)
after two years of preparation.
He was admitted together with Gustave Choquet, a~winner of
the  Concours G\'en\'eral in mathematics, and Marie-H\'el\`en, one of the
first females in the ENS.
The mathematical atmosphere of those years in the ENS
was determined by \'E. Borel, \'E.~Cartan, A.~Denjoy, M.~Fr\'echet, and
P.~Montel. The staff of the neighboring Coll\`ege de France included
H.~Lebesgue who delivered lectures and Hadamard who conducted  seminars.
It was in his student years that Schwartz had acquired the irretrievable and
permanent love to probability theory which grew from conversations with
his future father-in-law L\'evy.

After graduation from the ENS  Schwartz decided to be drafted in
the compulsory military service for two years. He had to stay in the army
in 1939--1940 in view of  the war times.
These years were especially hard for the young couple of the  Schwartzes.
It was unreasonable for Jews to stay in the occupied zone. The
Schwartzes had to escape from the native north and  manage
to survive on some modest financial support
that was offered in particular by Michelin, a world-renowned
tire company.
In 1941 Schwartz was in Toulouse for a short time and met H.~Cartan
and J.~Delsarte  who suggested that the young couple should move to
Clermont-Ferrand, the place of temporary residence of the group of professors
of Strasbourg University which had migrated from the German occupation.
These were J.~Dieudonn\'e, Ch.~Ehresmann, A.~Lichnerowicz, and
S.~Mandelbrojt. In Clermont-Ferrand  Schwartz completed his
Ph.~D. Thesis on approximation of a real function on
the axis by sums of exponentials.

Unfortunately, the war had intervened into
the mathematical fate of Schwartz. His family had
to change places under false  identities.
Curiously, in the time of the invention of
distributions in November of 1944
Schwartz used the identity of Selimartin.
The basics of Schwartz's theory appeared
in the {\it Annales} of the University of Grenoble in 1945.\footnote{Cp.~\cite{6}.}
Schwartz described the process of invention as
``cerebral percolation.''
After a year's stay in Grenoble, Schwartz acquired a position
in Nansy, plunging to the center of
``Bourbakism.'' It is well known that N.~Bourbaki resided
in Nancago, a~mixture of Nancy and Chicago.
A.~Weyl lived in Chicago, while Delsarte and Dieudonn\`e
were in Nancy.
Before long ~Schwartz  was enrolled in the group of Bourbaki.
In~1950 he was awarded with the Fields medal for
distribution theory. His now-celebrated two-volume set {\it Th\'eorie des Distributiones\/}
was printed a short time later.

In~1952 Schwartz returned back to Paris and began  lecturing in
Sorbonne;  and since 1959, in  the \'Ecole Politechnique
in company with his father-in-law L\'evy.
Many celebrities were the direct students of~Schwartz.
Among them we list A.~Grothendi\-eck, \hbox{J.-L.} Lions,
B.~Malgrange, and A.~Martineau.

Schwartz wrote:
``To discover something in mathematics is to overcome an inhibition
and a tradition. You cannot move forward if you are not subversive.''
This statement is in good agreement with
the very active and versatile public life of
Schwartz. He joined Trotskists in his green years, protesting
against the monstrosities of capitalism and Stalin's terror of the
1930s. Since then he had never agreed with anything that he viewed as
violation of human rights, oppression, or injustice.
He was very active in struggling against the
American war in Vietnam and the Soviet invasion in Afghanistan.
He fought for liberation of a few mathematicians that were
persecuted for political reasons, among them
Jose Luis Massera, Vaclav Benda, et al.

Schwartz was an outstanding lepidopterist and had collected more
than 20~000 butterflies. It is not by chance that the butterflies
are depicted on the soft covers of the second edition of
his  {\it Th\'eorie des Distributiones}.

Schwartz passed away in Paris on June~4, 2002.

\section{Advances of Distribution Theory}
Distribution theory stems from the intention to apply the technologies
of functional analysis to studying partial differential equations.
Functional analysis  rests on algebraization, geometrization, and socialization
of analytical  problems.
By socialization we usually mean the inclusion of a particular problem in
an appropriate class of its congeners.
Socialization enables us to erase the
``random features,''\footnote{This is a clich\'e with a century-old
history. The famous Russian symbolist Alexander
Blok (1880--1921) used the concept of random feature in his incomplete
poem ``Revenge'' as of 1910~\cite[p.~482]{Block}. The prologue of this poem contains the lines
that are roughly rendered in English as follows:\hfill\break
{\tiny
\obeylines\obeyspaces
\phantom{\quad\quad}You share the  gift of prudent measure\hfill\break
\phantom{\quad\quad}For what keen vision might perceive.\hfill\break
\phantom{\quad\quad}Erasing random features, treasure\hfill\break
\phantom{\quad\quad}The world of beauty to receive.}}
eliminating the  difficulties of the insurmountable specifics of a
particular problem. In the early 1930s
the merits of functional analysis were already demonstrated
in the area of integral equations. The time was ripe for the differential
equations  to be placed on  the agenda.

It is worth observing that the contemplations about
the nature of integration and differentiation underlie
most of the theories of the present-day functional analysis.
This is no wonder at all in view of
the key roles of these remarkable linear operations.
Everyone knows that integration possesses a few more attractive
features than differentiation: the integral is monotone and
raises smoothness. Derivation  lacks these nice properties completely.
Everyone knows as well that the classical
derivative yields a closed yet unbounded operator
(with respect to the natural uniform convergence topology that is
induced by the Chebyshev $\sup$-norm). The series of smooth functions
cannot be differentiated termwise in general, which diminishes the
scope of applications of analysis to differential equations.

There is practically no denying today that the concept of generalized
derivative occupies a central place in distribution theory.
Derivation is now treated as the operator that acts on the
nonsmooth functions according to the same integral laws  as the
procedure of taking the classical derivative.
It is exactly this approach that was pursued steadily
by Sobolev. The new turnpike led to
the stock of previously impossible differentiation formulas.
It turned out that each distribution possesses derivatives of all orders,
every series of distributions may be differentiated termwise
however often, and  many ``traditionally divergent''
Fourier series admit presentations by explicit  formulas.
Mathematics has acquired additional fantastic degrees of freedom,
which makes immortal the name of Sobolev as  a~pioneer of the
calculus of the twentieth century.

The detailed expositions of the new theory by Sobolev and Schwartz
had appeared practically at the same
time. In~1950 the first volume
of {\it Th\'eorie des Distributiones} was published in~Paris, while
Sobolev's book {\it Applications of Functional Analysis in
Mathematical Physics\/} was printed in
Leningrad. In~1962 the Siberian Division of the Academy of Sciences of the USSR
reprinted the book, while in 1963 it was translated into English by
the American Mathematical Society.
The second edition of the Schwartz book was published in 1966,
slightly enriched with a generalized version of the de Rham currents.
Curiously, Schwartz left the historical overview practically the same as
in the introduction to the first edition.

The new methods of distribution theory turned out so powerful as
to enable mathematicians to write down, in  explicit form,
the general solution of an arbitrary partial differential equation
with constant coefficients.
In fact, everything reduces to existence of fundamental solutions;
i.~e., to the case of the Dirac delta-function on the right-hand side
of the equation under consideration.
The existence of these solutions was already established
in~1953 and 1954 by B.~Malgrange and L.~Ehrenpreis independently
of each other.  However, it was only in~1994
that some formula for a fundamental solution was exhibited by
H.~K\"onig. Somewhat later
N.~Ortner and P.~Wagner found a more elementary formula.
Their main result is 
as follows:\footnote{Cp.~\cite{10} and~\cite[Theorem~2.3]{11}.}

\font\eightbf=cmbx8

\begin{itemize}
\item[]{\small\sl\indent
{\bf Theorem.} Assume that $P(\partial)\in\mathbb C[\partial]$,
where $P$ is a~polynomial of degree $m$.  Assume further that
$\eta\in\mathbb R^n$ and $P_m(\eta)\ne0$, where $P_m$ is the principal
part of~$P$; i.~e., $P_m=\sum\nolimits_{|\alpha|=m}a_{\alpha}\partial^{\alpha}$.
Then the distribution ~$E$ given as
$$
E:={\frac{1}{\overline{P_m(\eta)}}}
\int\limits_{\mathbb T}\lambda^m e^{\lambda\eta x}
{\mathfrak F}^{-1}_{\xi\to x}
\left(\frac{{\overline{P(i\xi+\lambda\eta)}}}
{{P(i\xi+\lambda\eta)}}\right)
{\frac{d\lambda}{2\pi i\lambda}}
$$
\noindent
is a fundamental solution of the operator $P(\partial)$.
Moreover, $E/\cosh(\eta x)\in\mathscr S'(\mathbb R^n)$.
}
\end{itemize}

\noindent
It stands to reason to inspect the structure of the formula
which reveals the role of the distributional Fourier transform
 $\mathfrak F$ and the Schwartz space $\mathscr S'(\mathbb R^n)$
 comprising tempered distributions.\footnote{Also known as
 ``generalized functions of slow growth.''}

The existence of a fundamental solution
of an arbitrary partial differential equation with constant coefficients
is reverently called the {\it Malgrange--Ehrenpreis Theorem.}
It is hard to overestimate this splendid achievement which
remains one of the splendid triumphs of the abstract theory of
topological vector spaces.

The road from solutions in distributions to standard solutions
lies through  Sobo\-lev spaces.
Study of the embeddings and traces of  Sobolev spaces
has become one of the main sections of the modern theory
of real functions.
Suffice it to mention S.~M.~Nikol$'$ski\u\i,
O.~V.~Besov, G.~Weiss, V.~P. Il$'$in, and V.~G.~Mazya in order
to conceive the greatness of this area of mathematical research.
The titles of dozen books mention Sobolev spaces, which
is far from typical in the present-day science.

The broad stratum of modern studies deals with applications of
distributions in mathematical and theoretical physics,
complex analysis, the theory of pseudodifferential operators,
Tauberian theorems, and other sections of mathematics.

The physical sources of distribution theory,
as well as the ties of the latter with theoretical physics,
are the topics of paramount importance. They require
a special scrutiny that falls beyond the scope of this
article.\footnote{Some historical details are collected
in~\cite{24}. Also see~\cite{25}. J.-M. Kantor
kindly made his article available to the author before publication
with a courteous cooperation of Ch.~Davis, Editor-in-Chief of
{\it The Mathematical Intelligencer}. It was the proposal of
Ch.~Davis that the article by J.-M.~Kantor be supplemented
with the short comments~\cite{26} and~\cite{27}.}
We will confine exposition to the concise historical comments
by~V.~S. Vladimirov:\footnote{Cited from the handwritten review for
the {\it Herald of the Russian Academy of Sciences}, dated
 as of December~10, 2003.}

\begin{itemize}
\item[]{\small\indent
It was already the creators of this theory,
S. L. Sobolev~\cite{4} and L.~Schwartz~\cite{17}
who studied the applications of the theory of generalized functions
in mathematical physics. After a conversation with
S.~L. Sobolev about generalized functions, N.~N. Bogolyubov
used the Sobolev classes~\cite{2} of test and generalized functions
$C^m_{\comp}$ and $(C^m_{\comp})^\ast$
in constructing his axiomatic quantum field theory
\cite{18}--\cite{20}. The same related to the  Wightman
axiomatics~\cite{21}.
Moreover, it is impossible in principle to
construct any axiomatics of quantum field theory without
generalized functions. Furthermore, in the theory
of the dispersion relations~\cite{22} that are
derived from the Bogolyubov axiomatics,
the generalized functions, as well as their generalizations--hyperfunctions,
appear as the boundary values of holomorphic functions of
(many) complex variables. This fact, together with
the relevant  aspects such as  Bogolyubov's ``Edge-of-the-Wedge''
 Theorem, essentially enriches the theory of generalized functions.}
\end{itemize}

\section{Various Opinions About the History of Distributions}

J.~Leray was one of the most prominent French mathematicians of
the twentieth century. He was awarded with the Lomonosov
Gold Medal together with Sobolev in~1988.
Reviewing the contributions
of Sobolev from~1930 to 1955 in the course of Sobolev's election
to the Academy of Sciences of the Institute of France in 1967,
J.~Leray wrote:\footnote{Cp.~\cite{12}.}

\begin{itemize}
\item[]{\small\indent
Distribution theory is now well developed
due to the theory of topological vector spaces and their duality
as well as the concept of tempered distribution which is one of the
important achievements of L.~Schwartz (Paris) which enabled him to
construct the beautiful theory of the Fourier transform for
distributions; G.~de Rham supplied the concept of distribution
with that of current which comprises the concepts of
differential form and topological chain;
L.~H\"ormander (Lund, Princeton), B.~Malgrange (Paris),
J.-L.~Lions (Paris) used the theory of distributions
to renew the theory of partial differential equations; while
P.~Lelong (Paris) established one of the fundamental properties
of analytic sets. The comprehensive two-volume treatise
by  L.~Schwartz and even more comprehensive five-volume
treatise\footnote{In fact, the series consists
of 6 volumes.}
 by Gelfand and Shilov (Moscow) are the achievements
of so great an importance that even the French contribution
deserves the highest awards of our community. The applications
of distribution theory in all areas of mathematics, theoretical
physics, and numerical analysis remind of the dense forest
hiding the tree whose seeds  it has grown from.
However, we know that if Sobolev had fail to make his discovery
about 1935 in Russia, it would be committed in France by 1950
and somewhat later in Poland; the USA also flatters itself
that they would make this discovery in the same years:
The science and art of mathematics would be late only by 15 years as
compared with Russia\dots.}
\end{itemize}

\noindent
In sharp contrast with this appraisal, we cite
F.~Tr\'eves who wrote in the memorial article about
Schwartz in October 2003 
as follows:\footnote{Cp.~\cite[p.~1076]{13}.}

\begin{itemize}
\item[]{\small\indent
The closest any mathematician of the 1930s ever
came to the general definition of a distribution is
Sobolev in his articles [Sobolev, 1936] and [Sobolev,
1938]\footnote{These are references to the articles in {\it Sbornik} [2, 3].}
 (Leray used to refer to ``distributions, invented
by my friend Sobolev''). As a matter of fact,
Sobolev truly defines the distributions of a given,
but arbitrary, finite order $m$: as the {\it continuous
linear functionals} on the space $C^m_{\comp}$ of
compactly supported functions of class $C^m$. He
keeps the integer $m$ fixed; he never considers the
intersection $C^{\infty}_{\comp}$  of the spaces $C^m_{\comp}$
for all $m$. This
is all the more surprising, since he proves that
$C^{m+1}_{\comp}$C  is dense in $C^m_{\comp}$     by the Wiener procedure of
convolving functions $f\in C^m_{\comp}$ with a sequence of
functions belonging to $C^{\infty}_{\comp}$! In connection with
this apparent blindness to the possible role of
mentioned to Henri Cartan his inclination to use
the elements of $C^{\infty}_{\comp}$ as test functions, Cartan
tried to dissuade him: ``They are too freakish ({\it trop
monstrueuses}).''}

\item[]{\small\indent
Using transposition, Sobolev defines the multiplication
of the functionals belonging 	to $C^m_{\comp}$
by the functions belonging to $C^m$ and the differentiation
of those functionals: $d/dx$ maps $(C^m_{\comp})^\ast$
into $(C^{m+1}_{\comp})^\ast$. But again there is no mention of
Dirac $\delta(x)$ nor of convolution, and no link is made
with the Fourier transform. He limits himself to
applying his new approach to reformulating and
solving the Cauchy problem for linear hyperbolic
equations. And he does not try to build on his
remarkable discoveries. Only after the war does
he invent the Sobolev spaces $H^m$ and then only for
integers $m\ge 0$. Needless to say, Schwartz had not
read Sobolev's articles, what with military service
and a world war (and Western mathematicians'
ignorance of the works of their Soviet colleagues).
There is no doubt that knowing those articles would
have spared him months of anxious uncertainty.}
\end{itemize}

\noindent
F.~Tr\'eves should be honored for drifting aside
from the practice of evaluating publications from
what they lack when he wrote somewhat later about that which
made the name of Schwartz immortal:\footnote{Ibid., p.~1077.}

\begin{itemize}
\item[]{\small\indent
 Granted that Schwartz might have been
replaceable as the inventor of distributions, what
can still be regarded as his greatest contributions
to their theory? This writer can mention at least two
that will endure: (1) deciding that the Schwartz
space $\mathscr S$  of rapidly decaying functions at infinity and
its dual $\mathscr S'$    are the ``right'' framework for Fourier
analysis, (2) the Schwartz kernel theorem.}
\end{itemize}

\noindent
The Tr\'eves opinion coincides practically verbatim with
the narration of Schwartz in his autobiography
published firstly in~1997.
Moreover, Schwartz had even remarked there about Sobolev
that\footnote{Cp.~\cite[p.~222]{9}.}

\noindent
\begin{itemize}
\item[]{\small\indent
he did not develop his theory
in view of general applications, but with a precise goal:
he wanted to define the generalized solution of a partial
differential equation with a second term and initial conditions.
He includes the initial conditions in the second term in the
form of functionals on the boundary and obtains in this way
a remarkable theorem on  second order hyperbolic partial
differential equations. Even today this remains one of the
most beautiful applications of the theory of distributions,
and he  found it in a~rigorous manner. The astounding thing is
that he stopped at this point. His 1936 article, written in French,
is entitled ``Nouvelle m\'ethode \`a r\'esoudre de probl\`eme
de Cauchy pour les \'equations lin\'eares hyperboliques normales.''
After this article, he did nothing further in this fertile
direction. In other words, Sobolev himself did not fully
understand the importance of his discovery.}
\end{itemize}

\noindent
It is impossible to agree with these opinions.
Rather strange is to read about the absence of any mention of
the Dirac delta-function among the generalized functions of Sobolev,
since $\delta$ obviously belongs to each of
the spaces~$(C^m_{\comp})^\ast$.

Disappointing is the total neglect of the classical treatise
of Sobolev~\cite{5} which was a deskbook of many specialists in functional analysis
and partial differential equations
for decades.\footnote{Published in~1950 by Leningrad State University,
reprinted in~1962 by the Siberian Division of the Academy of Sciences of
the USSR in~Novosibirsk, and translated into English by the American Mathematical Society
in~1963. The third Russian edition was printed by the Nauka Publishers in~1988.}
Finally, Schwartz was not recruited in 1997 and did not participate
in a`world war. Therefore, there were some other reasons
for him to neglect the Sobolev book~\cite{6}
which contains  the principally new
Sobolev based his pioneering results in numerical integration
on developing the theory of the Fourier transform of
distributions which was created by~Schwartz.

Prudent in the appraisals,  exceptionally tactful, and modest in his ripe years,
Sobolev always abstained from any bit of details of the history of distribution
theory neither in private conversations nor in his numerous writings.
The opinion that he decided worthy to be left to the future
generations about this matter transpires in his concise comments on
the origins of distribution theory in his book~\cite[Ch.~8]{6}
which was printed in~1974:

\begin{itemize}
\item[]{\small\indent
\quad\quad The generalized functions are ``ideal elements''
that complete the classical function spaces
in much the same way as the real numbers complete
the set of rationals.
}
\item[]{\small\indent
\quad\quad In this chapter we concisely present the theory
of these functions which we need in the sequel.
We will follow the way of presentation
close to that which was firstly used by the author
in~1935 in~[16].\footnote{This is a~reference to
the article of~1936 in~{\it Sbornik}~\cite{2}.}
The theory of generalized functions
was further developed by L.~Schwartz~[21]
who has  in particular considered and studied the Fourier transform
of a~generalized function.\footnote{Cp.~[25, p.~355].
This is a curious misprint: the correct reference to Schwartz's  two-volume set
should be~[47].}
}
\item[]{\small\indent
\quad\quad Historically, the generalized function
had appeared explicitly in the studies in theoretical physics
as well as in the works of J.~Hadamard, M.~Riesz, S.~Bochner,
et al.
}
\end{itemize}

\noindent
Therefore, we can agree only in part with the
following statement by~Schwartz~\cite[p.~236]{9}: 
\begin{itemize}
\item[]{\small\indent
Sobolev and I and all the others who came before us
were influenced by out time, our environment and our own previous work.
This makes it less glorious, but since we were both ignorant
of the work of many other people, we still had to develop plenty of
originality.}
\end{itemize}

\noindent
Most mathematicians  agree that Israel Gelfand could be ranked as
the best arbiter in distribution theory.
The  series {\it Generalized Functions} written by him and his
students was started in the mid 1950s and remains one of the heights
of the world mathematical literature, the encyclopedia of distribution
theory. In the preface to the first edition of the first volume of this series,
Gelfand wrote:\footnote{Cp.~\cite{14}.}

\begin{itemize}
\item[]{\small\indent
It was S.~L.~Sobolev who introduced generalized functions
in explicit and now generally accepted form
in~1936\dots. The monograph of Schwartz {\it Th\'eorie des Distributiones}
appeared in~1950--1951. In  this book  Schwartz
systemized the theory of generalized functions,
interconnected all previous approaches, laid  the theory
of topological linear spaces in the foundations of the theory
of generalized functions, and obtain a number of essential and
far-reaching results. After the publication of
{\it Th\'eorie des Distributiones},  the generalized functions
won exceptionally swift and wide popularity just in two or three years.}
\end{itemize}

\noindent
This is an accurate and just statement. We may agree with it.

\section{Classicism and Romanticism}

Pondering over the fates of~Sobolev and~Schwart, it is impossible
to obviate the problem of polarization of the opinions about
the mathematical discovery of these scholars.
The hope is naive  that this problem will ever received
a simple and definitive answer that satisfies and convinces everyone.
It suffices to consider the available experience that concerns
other famous pairs of mathematicians whose fates and contributions
raise the quandaries that sometimes lasted for centuries and resulted
in the fierce clashes of opinions up to the present day.
The sources of these phenomena seem of a rather universal provenance
that is not concealed in the particular  personalities
but resides most probably in the nature of mathematical creativity.

Using quite a risky analogy, we may say that mathematics
has some features  associated with the trends of artistry
which  are referred to traditionally as classicism and romanticism.
It is hard to fail discerning the classic lineaments
of the Hellenistic tradition in the writings of
Euclid, Newton, Bolyai, Hilbert, and Bourbaki.
It is impossible to fail to respond to the accords of the
romantic anthem of the human genius which sound
in the pages of the writings of Diophant, Leibniz, Lobachevski\u\i,
Poincar\'e, and Arnold.

The magnificent examples of  mathematical classicism and romanticism
glare from the creative contributions of
Sobolev and~Schwartz.  These giants and their achievements
will remain with us for ever.

\bibliographystyle{plain}

\end{document}